\documentclass[12pt]{article} 
\usepackage{amsfonts,latexsym}

\newcounter{environment}[section]
\renewcommand{\theenvironment}{%
\arabic{section}.\arabic{environment}}

{\begin{rm}\refstepcounter{environment}{\textbf\theenvironment\
\bf Definition.~~}}%
{\end{rm}}

{\begin{rm}\refstepcounter{environment}{\textbf\theenvironment\
\bf Example.~~}}%
{\end{rm}}

\newenvironment{conjecture}%
{\begin{rm}\refstepcounter{environment}{\textbf\theenvironment\
\bf Conjecture.~~}}%
{\end{rm}}

\newenvironment{proposition}%
{\begin{rm}\refstepcounter{environment}{\textbf\theenvironment\
\bf Proposition.~~}}%
{\end{rm}}

{\begin{rm}\refstepcounter{environment}{\textbf\theenvironment\
\bf Proposition and Definition.~~}}%
{\end{rm}}

\newenvironment{theorem}%
{\begin{rm}\refstepcounter{environment}{\textbf\theenvironment\
\bf Theorem.~~}}%
{\end{rm}}

\newenvironment{corollary}%
{\begin{rm}\refstepcounter{environment}{\textbf\theenvironment\
\bf Corollary.~~}}%
{\end{rm}}

{\begin{rm}\refstepcounter{environment}{\textbf\theenvironment\
\bf Lemma.~~}}%
{\end{rm}}

\begin{document}
\newcommand{\qbc}[2]{{\left [{#1 \atop #2} \right ]}}
\newcommand{\beq}{\begin{equation}}
\newcommand{\eeq}{\end{equation}}
\newcommand{\fq}{\mathbb{F}_q}
\newcommand{\ale}{\mbox{$<\!\!\! *$}}
\newcommand{\aleq}{\mbox{$\leq\!\! *$}}
\newcommand{\pbar}{\bar{\mathbb{P}}}
\newcommand{\be}{\begin{enumerate}}
\newcommand{\ee}{\end{enumerate}}
\newcommand{\bea}{\begin{eqnarray}}
\newcommand{\eea}{\end{eqnarray}}
\newcommand{\beas}{\begin{eqnarray*}}
\newcommand{\eeas}{\end{eqnarray*}}
\newcommand{\sym}{\mathrm{Sym}}
\newcommand{\gln}{\mathrm{GL}(n,q)}
\newcommand{\ds}{\displaystyle}

\begin{centering}
{\huge\bf Spanning Trees and a Conjecture of Kontsevich}\\[.2in]
Richard P. Stanley\footnote{Partially supported by NSF grant
\#DMS-9743966.}\\ 
Department of Mathematics\\
Massachusetts Institute of Technology\\
Cambridge, MA 02139\\
\emph{e-mail:} rstan@math.mit.edu\\[.2in]
version of 7 November 1998\\[.2in]
\textbf{Running title}: Spanning Trees and a Conjecture of
Kontsevich\\[.3in] 

\end{centering}
\vskip 10pt
\section{Introduction.}
Let $G$ be a (finite) graph with vertex set
$V(G)=\{v_1,\dots,v_n\}$ and edge set $E(G)=\{e_1,\dots,e_s\}$. For
now we allow loops (edges from a vertex to itself) and multiple edges.
For each edge $e$ of $G$ associate an indeterminate $x_e$. If
$S\subseteq E(G)$ then let
  $$ x^S = \prod_{e\in S} x_e. $$
Let $T$ be a spanning tree of $G$, and let $\bar{T}$ denote the (edge)
complement of $T$. Define polynomials $P_G(x)$ and $Q_G(x)$ in the
variables $x=(x_e)_{e\in E(G)}$ by
  \bea P_G(x) & = & \sum_T x^{\bar{T}}\nonumber \\
        Q_G(x) & = & \sum_T x^T, \label{eq:qgdef} \eea
where both sums range over all spanning trees $T$ of $G$. For
instance, if $G$ is a four-cycle with edge set $\{1,2,3,4\}$, then
  \beas P_G(x) & = & x_1+x_2+x_3+x_4\\
        Q_G(x) & = & x_1x_2x_3+x_1x_2x_4+x_1x_3x_4+x_2x_3x_4. \eeas
Note that if $G$ isn't connected, then $P_G(x)=Q_G(x)=0$.
Note also that for any $G$, we have
  \beq Q_G(x) = \left( \prod_{e\in E(G)} x_e\right)P_G(1/x), 
    \label{eq:1} \eeq
where $P_G(1/x)$ denotes the result of substituting $1/x_e$ for $x_e$
in $P_G(x)$ (for all $e$).

Let $q$ be a prime power, and let $f_G(q)$ (respectively, $g_G(q)$)
denote the number of ways of substituting elements of the finite field
$\mathbb{F}_q$ for the variables $x_e$ so that $P_G(x)\neq 0$
(respectively, $Q_G(x)\neq 0$). For instance, if $G$ is a four-cycle as
above and we want $P_G(x)\neq 0$, then $x_1,x_2,x_3$ can be arbitrary,
and then there are $q-1$ choices for $x_4$. Hence $f_G(q)=q^3(q-1)$. 
If we want $Q_G(x)\neq 0$ then it takes a little more work to see that
$g_G(q)=q(q-1)(q^2-2)$. More generally, if $C_n$ denotes an $n$-cycle
then it is not hard to check (as was done first in the case of
$g_{C_n}(q)$ by J. Stembridge) that
  \beas f_{C_n}(q) & = & q^{n-1}(q-1)\\[.05in]
        g_{C_n}(q) & = & n(q-1)^{n-1} +
     (q-1)^n - (q-1)^{n-1} + (q-1)^{n-2}\\
      & & \qquad+\cdots+ (-1)^{n-1}(q-1).
  \eeas
Note that if $G$ isn't connected then $f_G(q)=g_G(q)=0$, since
$P_G(x)=Q_G(x)=0$. 

In a lecture delivered at the Rutgers
University Gelfand Seminar on December 8, 1997, M. Kontsevich stated
the conjecture, in connection with the evaluation of certain
integrals appearing in perturbative quantum field theory, that $f_G(q)$ 
is a ``universal polynomial'' in $q$, i.e., a polynomial in $q$
independent of the characteristic $p$ of the field $F_q$. We
have been unable to resolve Kontsevich's conjecture, but in
Sections~\ref{sec:cg} and \ref{sec:cgg} we present evidence in its
favor while in Section~\ref{sec:neg} we present evidence against
it. Let us mention that John Stembridge has verified that $g_G(q)$ is
a polynomial for all but two graphs with twelve edges. Moreover, for
these two exceptional graphs he has computed enough values of $g_G(q)$
to overdetermine the polynomial, without contradicting that $g_G(q)$
is a polynomial. Hence it is quite likely that $g_G(q)$ is also a
polynomial for these two graphs.


\section{Some general observations.}
Given the graph $G$, let $S$ be a subset of $E(G)$. Define
$f_{G,S}(q)$ (respectively, $f_{G,S}^+(q)$) to be the number of ways
of substituting elements of $\fq$ for the variables $x_e$ such that
$x_e=0$ if $e\in S$ (respectively, if and only if $e\in S$) so that
$P_G(x) \neq 0$. Similarly define $g_{G,S}(q)$ and
$g_{G,S}^+(q)$ using $Q_G(x)$. In particular,
$f_{G,\emptyset}(q)=f_G(q)$ and $g_{G,\emptyset}(q)=g_G(q)$. Now 
  \beas f_{G,S}(q) & = & \sum_{T\supseteq S} f_{G,T}^+(q)\\
        g_{G,S}(q) & = & \sum_{T\supseteq S} g_{G,T}^+(q). \eeas
Hence by the Principle of Inclusion-Exclusion we have
  \beas f_{G,S}^+(q) & = & \sum_{T\supseteq S}
     (-1)^{\#(T-S)}f_{G,T}(q)\\ 
        g_{G,S}^+(q) & = & \sum_{T\supseteq S} (-1)^{\#(T-S)}
        g_{G,T}(q). \eeas 
Now note that $g_{G,S}(q)=g_H(q)$, where $H$ is the spanning subgraph
of $G$ with edge set $E(G)-S$. Similarly, if $S$ is acyclic (contains
no cycle) then $f_{G,S}(q)=f_K(q)$ where $K$ denotes $G$ with the
edges in $S$ contracted to points. On the other hand, if $S$ contains
a cycle then $f_{G,S}(q)=0$. Finally observe from 
equation~(\ref{eq:1}) that for any graph $H$ we have
 $$ f_{H,\emptyset}^+(q) = g_{H,\emptyset}^+(q). $$
From these observations we obtain the following result. 

\begin{proposition}
\emph{Let $n\geq 1$. Then $f_G(q)$ is a universal
  polynomial in $q$ for all graphs $G$ on at most $n$ vertices if and
  only if the same is true for $g_G(q)$.}
\end{proposition}

From now on we will deal only with $Q_G(x)$ and $g_G(q)$. Note also
that if $G'$ denotes $G$ with either one new loop added to a vertex or
one edge replaced by two edges, then
  $$ g_{G'}(q) = q\cdot g_G(q). $$
Hence it suffices to assume from now on that $G$ is \emph{simple},
i.e., has no loops or multiple edges.

John Stembridge has pointed out that a minimal counterexample to
Kontsevich's conjecture must have edge-connectivity at least three,
i.e, the graph cannot be disconnected by the removal of two
edges. Indeed, suppose that there are two edges $e$ and $e'$ whose
removal disconnects $G$, leaving the disjoint union of $G_1$ and
$G_2$. Let $G'$ denote the graph obtained from $G$ by contracting $e$,
and let $G''$ be obtained by further contracting $e'$. Then it is easy
to check that
  $$ g_G(q) = q\cdot g_{G_1}(q)g_{G_2}(q)+(q-2)g_{G'}(q)+(q-1)g_{G''}(q). $$

As a final remark, recall that by the rationality
of the zeta function of an algebraic variety over a finite field (due
to Dwork \cite{dwork}), for a fixed prime power $q$ there exist
algebraic integers $\alpha_1,\dots, \alpha_r$ and
$\beta_1,\dots,\beta_s$ such that for all $m\geq 1$,
  $$ g_G(q^m) = \sum \alpha_i^m - \sum \beta_j^m. $$
Since the $\alpha_i$'s and $\beta_j$'s are algebraic integers, it is
easy to deduce the following consequence.

\begin{proposition} \emph{If $g_G(q)$ is a polynomial in $q$ with}
  rational \emph{coefficients, then in fact $g_G(q)$ has} integer
  \emph{coefficients}. 
\end{proposition}

\section{The Matrix-Tree Theorem and nonsingular symmetric matrices.}
The fundamental tool for our results is the Matrix-Tree Theorem. This
result was stated by J. J. Sylvester in 1857. The first proof was
published by C. W. Borchardt in 1860. The Matrix-Tree Theorem is often
attributed to A. Cayley because he cited Sylvester's work in 1856
before it was published. For an exposition of the Matrix-Tree Theorem
and more precise references, see Chapter~5.6 and the Notes
to Chapter 5 of \cite{ec2}. Let $G$ be a graph without loops or
multiple edges on the vertex set $\{v_1,\dots,v_n\}$, and as above
associate the indeterminate $x_e$ with every edge $e$ of $G$. Let $L
=L(G)=\left( L_{ij}\right)_1^n$ be the $n\times n$ matrix defined by
  $$ L_{ij} = \left\{ \begin{array}{rl} -x_e, & \mbox{if $i\neq j$
        and $e$ has vertices $v_i$ and $v_j$}\\ 
     0, & \mbox{if $i\neq j$ and there is no edge with vertices $v_i$
        and $v_j$}\\
     \sum_e x_e, & \mbox{if $i=j$ and $e$ ranges over all edges
        incident to $v_i$.} \end{array} \right. $$ 

Let $L_0=L_0(G)$ denote $L$ with the last row and column removed. We
call $L$ the \emph{generic Laplacian matrix} of $G$, and $L_0$ the
\emph{reduced} generic Laplacian matrix of $G$.

\begin{theorem} (the Matrix-Tree Theorem) \emph{We have}
  $$ Q_G(x) = \det L_0(G). $$
\end{theorem}

\begin{corollary} \label{cor:rlap} 
\emph{Let $L_0$ be the reduced generic Laplacian 
matrix of the connected graph $G$. Then $g_G(q)$ is the number of
solutions to $\det L_0\neq 0$ over $\mathbb{F}_q$.}
\end{corollary}

We say that a vertex $v$ of the simple graph $G$ is an \emph{apex} if
$v$ is incident to every other vertex of $G$. For graphs with apexes
Corollary~\ref{cor:rlap} has the following variant.

\begin{theorem} \label{thm:apex}
\emph{Let $G$ be a simple graph with vertices
$v_1,\dots,v_n$ such that $v_n$ is an apex. Then $g_G(q)$ is equal
to the number of $(n-1)\times(n-1)$ nonsingular symmetric matrices $M$
over $\mathbb{F}_q$ such that $M_{ij}=0$ whenever $i\neq j$ and $G$
has no edge between $v_i$ and $v_j$.}
\end{theorem}

\textbf{Proof.} Let $e_i$ be the unique edge of $G$ with vertices
$v_i$ and $v_n$, for $1\leq i\leq n-1$. The $(i,i)$-entry $(L_0)_{ii}$
of the reduced generic Laplacian matrix $L_0$ has the form $x_{e_i} +$
other terms, and nowhere else does $x_{e_i}$ appear in $L_0$. Hence we
can replace $(L_0)_{ii}$ with $x_{e_i}$ without affecting the set of
$q^{\#E}$ matrices we obtain from $L_0$ by letting the $x_e$'s assume
all possible values in $\mathbb{F}_q$. Similarly we don't affect this
set by changing the signs of the off-diagonal entries. But then $L_0$
becomes a symmetric matrix $L_0^*$ whose entries are generic except
that $(L_0^*)_{ij}=0$ whenever $i\neq j$ and $G$ has no edge between
$v_i$ and $v_j$, and the proof follows. $\ \Box$

\section{The complete graph.} \label{sec:cg}
Theorem~\ref{thm:apex} allows us to evaluate $g_G(q)$ explicitly for
certain graphs $G$. We first consider the complete graph $K_n$, with
$n$ vertices and one edge between every pair of distinct
vertices. Hence by Theorem~\ref{thm:apex}, $g_{K_n}(q)$ is just the
total number of nonsingular $(n-1)\times (n-1)$ symmetric matrices over
$\mathbb{F}_q$. This number was first computed for $q$ odd by L. Carlitz
\cite[Thm.\ 3]{carlitz} as part of a much more general result. A simpler
proof valid for any $q$ was later given by J. MacWilliams \cite[Thm.\
2]{macw}. We will sketch the 
proof of MacWilliams and a second proof based on orthogonal geometry
over $\mathbb{F}_q$, since both proofs will lead to generalizations.

\begin{theorem} \label{thm:kn}
\emph{We have}
  \beq g_{K_n}(q) = \left\{ \begin{array}{rl}  
   q^{m(m-1)}(q-1)(q^3-1)\cdots(q^{2m-1}-1),\quad n=2m\\[.1in]
   q^{m(m+1)}(q-1)(q^3-1)\cdots(q^{2m-1}-1),\quad n=2m+1.
  \end{array} \right. \label{eq:gknq} \eeq 
\end{theorem}

\textbf{First proof} (J. MacWilliams). It is more convenient to
consider $g_{K_{n+1}}(q)$ rather than $g_{K_n}(q)$, so that we are
enumerating $n\times n$ invertible symmetric matrices over
$\fq$.  Let $h(n,r)$ denote the number of $n\times n$
symmetric matrices $M$ over $\fq$ of rank $r$. We claim that
  \beq h(n,r) = \left\{ \begin{array}{rl} \ds\prod_{i=1}^s
     \frac{q^{2i}}{q^{2i}-1}\cdot\prod_{i=0}^{2s-1}\left(
     q^{n-i}-1\right), & 0\leq r=2s\leq n\\[.15in] \ds
    \prod_{i=1}^s
     \frac{q^{2i}}{q^{2i}-1}\cdot\prod_{i=0}^{2s}\left(
     q^{n-i}-1\right), & 0\leq r=2s+1\leq n. \end{array}
     \right. \label{eq:ranksym} \eeq
An $(n+1)\times(n+1)$ symmetric matrix may be written as
  $$ N = \left[ \begin{array}{cc} \beta & y\\ y^t & M
     \end{array} \right], $$
where $M$ is an $n\times n$ symmetric matrix, $\beta\in\fq$, and $y
\in\fq^n$. Elementary linear algebra arguments (given explicitly in
\cite{macw}) show that from a particular $M$ of rank $r$ we obtain:
  \begin{itemize} \item $q^{n+1}-q^{r+1}$ matrices $N$ of rank $r+2$,
   \item $(q-1)q^r$ matrices $N$ of rank $r+1$,
   \item $q^r$ matrices $N$ of rank $r$,
   \item no matrices of other ranks.
  \end{itemize}
There follows the recurrence
  \beq h(n+1,r) = q^rh(n,r)+(q-1)q^{r-1}h(n,r-1)+(q^{n+1}
     -q^{r-1})h(n,r-2). \label{eq:rec} \eeq
One can check that the solution to this recurrence satisfying the
initial conditions $h(n,0)=1$ and $h(n,r)=0$ for $r>n$ is given by
(\ref{eq:ranksym}). The proof follows from the case $r=n$. $\ \Box$
   
\textbf{Second proof}. \emph{Case 1}: $q$ odd. Let $\sym^+(n,q)$
(respectively, $\sym^-(n,q)$) denote the set of all $n\times n$
nonsingular symmetric matrices over $\mathbb{F}_q$ whose determinant
is a square (respectively, a nonsquare) in $\mathbb{F}_q$. Let
$\Omega^+(n,q)$ denote the group of all matrices $A\in\mathrm{GL}(n,q)$
satisfying $AA^t=I$. (We will be dealing with various groups closely
related to the orthogonal groups $O(n,q)$. We will use the notation
$\Omega$ rather than $O$ to make clear that our groups are related but
in general not equal to the usual orthogonal groups.) By standard
results concerning orthogonal geometry over a finite field
(implicit in \cite{dickson}), the map $f:\mathrm{GL}(n,q)\rightarrow
\sym^+(n,q)$ defined by $f(A)=AA^t$ is surjective, and $f^{-1}(AA^t)=
A\cdot \Omega^+(n,q)$, the left coset of $\Omega^+(n,q)$ in GL$(n,q)$
containing 
$A$. Hence all fibers $f^{-1}(B)$ have cardinality $\#\Omega^+(n,q)$, so
  $$ \#\sym^+(n,q) = \frac{\#\mathrm{GL}(n,q)}{\#\Omega^+(n,q)}. $$
Similarly, let $\alpha$ be a fixed nonsquare in $\mathbb{F}_q$, and
let $\Omega^-(n,q)$ denote the group of all matrices $A\in\mathrm{GL}(n,q)$
satisfying $ADA^t=D$, where
$D=\mathrm{diag}(\alpha,1,1,\dots,1)$. Then the map
$f:\mathrm{GL}(n,q)\rightarrow 
\sym^-(n,q)$ defined by $f(A)=ADA^t$ is surjective, and $f^{-1}(ADA^t)=
A\cdot \Omega^-(n,q)$. Hence all fibers $f^{-1}(B)$ have cardinality
$\#\Omega^-(n,q)$, so 
  $$ \#\sym^-(n,q) = \frac{\#\mathrm{GL}(n,q)}{\#\Omega^-(n,q)}. $$
Since $\sym(n,q)=\sym^+(n,q)\bigcup\sym^-(n,q)$, there follows
  \beq \#\sym(n,q) = \frac{\#\mathrm{GL}(n,q)}{\#\Omega^+(n,q)} +
     \frac{\#\mathrm{GL}(n,q)}{\#\Omega^-(n,q)}. \label{eq:symnq} \eeq
The order of GL$(n,q)$ is well-known and easily seen to be 
  $$ \#\mathrm{GL}(n,q) = (q^n-1)(q^n-q)\cdots(q^{n-1}-1). $$
Moreover, the orders of $\Omega^+(n,q)$ and $\Omega^-(n,q)$ were computed by
Dickson \cite[Thm.\ 172]{dickson}, as follows. 

  \beas \#\Omega^+(n,q) & = & \left\{ \begin{array}{rl} \ds
    2q^{m(m-1)}(q^m-1)\prod_{i=1}^{m-1} (q^{2i}-1), & n=2m,\ 
        q\equiv 1\ (\mathrm{mod}\ 4)\\[.05in] \ds
    2q^{m(m-1)}(q^m-(-1)^m)\prod_{i=1}^{m-1} (q^{2i}-1), & n=2m,\ 
        q\equiv 3\ (\mathrm{mod}\ 4)\\[.05in]
    \ds 2q^{m^2}\prod_{i=1}^m (q^{2i}-1), & n=2m+1. 
   \end{array} \right. \\[.15in]
 \#\Omega^-(n,q) & = & \left\{ \begin{array}{rl} \ds
    2q^{m(m-1)}(q^m+1)\prod_{i=1}^{m-1} (q^{2i}-1), & n=2m,\ 
        q\equiv 1\ (\mathrm{mod}\ 4)\\[.05in] \ds
    2q^{m(m-1)}(q^m+(-1)^m)\prod_{i=1}^{m-1} (q^{2i}-1), & n=2m,\ 
        q\equiv 3\ (\mathrm{mod}\ 4)\\[.05in]
    \ds 2q^{m^2}\prod_{i=1}^m (q^{2i}-1), & n=2m+1. 
   \end{array} \right. \eeas

Substituting these numbers into equation (\ref{eq:symnq}) (after
replacing $n$ by $n-1$) yields (\ref{eq:gknq}) (when $q$ is odd).

\emph{Case 2}: $q$ even. This case is analogous to the odd case, but
the details are somewhat different. When $n$ is odd, it follows from
\cite[Thm.\ 7]{albert} that the map
$f:\mathrm{GL}(n,q)\rightarrow \sym(n,q)$ defined by $f(A)=AA^t$ is
surjective, with $\#f^{-1}(AA^t)=A\cdot \Omega(n,q)$ (where $\Omega(n,q)=
f^{-1}(I)$). Hence
  $$ \#\sym(n,q) = \frac{\#\mathrm{GL}(n,q)}{\#\Omega(n,q)}. $$
Dickson \cite[p.\ 206]{dickson} showed that
  $$ \#\Omega(2m+1,q) = q^{m^2} \prod_{i=1}^m (q^{2i}-1), $$
so (\ref{eq:gknq}) follows in this case.

When $n=2m$, let $\sym^+(n,q)$ (respectively, $\sym^-(n,q)$) denote
the set of $n\times n$ nonsingular matrices over $\fq$ with at least
one nonzero entry on the main diagonal (respectively, with all 0's on
the main diagonal). (When $n$ is odd we have $\sym^-(N,q)=\emptyset$,
since a symmetric matrix of odd order with zero diagonal over a field
of characteristic two is singular.)  It was shown by Albert
\cite[Thm.\ 7]{albert} that the map 
$f:\gln\rightarrow\sym^+(n,q)$ defined by $f(A)=AA^t$ is
surjective. Let $E$ be the direct sum of $m$ copies of the matrix
$\left[ \begin{array}{cc} 0 & 1\\ 1 & 0\end{array} \right]$. The map
$f:\gln\rightarrow \sym^-(n,q)$ defined by $f(A)=AEA^t$ is
surjective. If $\Omega^-(n,q)=f^{-1}(I)$, then $f^{-1}(AEA^t) = A\cdot
\Omega^-(n,q)$. Hence reasoning as before gives
  $$ \#\sym(n,q) = \frac{\#\mathrm{GL}(n,q)}{\#\Omega^+(n,q)} +
     \frac{\#\mathrm{GL}(n,q)}{\#\Omega^-(n,q)}. $$
It follows from the work of Dickson \cite[Ch.\ VIII]{dickson} that 
  \beas \#\Omega^+(2m,q) & = & q^{m^2}\prod_{i=1}^{m-1} (q^{2i}-1)\\[.1in]
        \#\Omega^-(2m,q) & = & q^{m^2}\prod_{i=1}^m (q^{2i}-1), \eeas 
from which we obtain (\ref{eq:gknq}) in this final case. $\ \Box$

Note that the first proof of Theorem~\ref{thm:kn} makes it clear from
the start that $g_{K_n}(q)$ is a polynomial, while in the second proof
(especially when $n$ is even) it appears somewhat miraculous that the
computations in odd and even characteristics lead to the same final
answer. The fact that the two cases yield the same answer boils down
to the identity
  $$ \frac 12\left( \frac{1}{q^m-1}+\frac{1}{q^m+1}\right) =
      \frac{1}{q^m}+\frac{1}{q^m(q^{2m}-1)}. $$

\section{Some generalizations of the complete graph.} \label{sec:cgg}
The two proofs we gave of Theorem~\ref{thm:kn} can be extended to more
general results. For the first generalization, let $G$ be an
$n$-vertex graph (without loops or multiple edges). Let $L_0$ denote
the reduced generic Laplacian matrix of $G$, with respect to some
vertex $v$ indexing the last row and column of $L$.  Write $h(G,r)$
for the number of ways of evaluating $L_0$ over $\fq$ (i.e., the
number of ways to substitute elements of $\fq$ for the variables
appearing in $L_0$) such that a matrix of rank $r$ is obtained. Thus
if $v$ is an apex of $G$, then by Theorem~\ref{thm:apex} we have
$h(G,n-1)=g_G(q)$.

\begin{theorem} \label{thm:genmcw}
\emph{Let $G$ be an $n$-vertex graph with an apex, and let $G^*$
  denote $G$ with an apex adjoined (so $G^*$ has at least two
  apexes). Then} 
  \beq  h(G^*,r) = q^rh(G,r)+(q-1)q^{r-1}h(G,r-1)+(q^n-q^{r-1})
      h(G,r-2). \label{eq:gast} \eeq
\end{theorem}

\textbf{Proof.} The proof is essentially the same as the first proof
of Theorem~\ref{thm:kn}. Let the vertices of $G^*$ be $v_1, \dots,
v_{n+1}$, where $v_1$ and $v_{n+1}$ are apexes. Let $e_i$ denote the
edge from $v_1$ to $v_i$ for $2\leq i\leq n+1$, and write $x_i$ for
$x_{e_i}$. We then have
  $$ L_0(G^*) = \left[ \begin{array}{cc} \beta & y\\ y^t & L_0(G)
     \end{array} \right], $$
where $y=-(x_2,\dots,x_n)$ and $\beta=x_2+\cdots+x_n+x_{n+1}$. Since
$x_{n+1}$ appears in $L_0(G^*)$ only in the entry $\beta$, we may
replace $\beta$ by $x_{n+1}$ without affecting the set of matrices we
obtain from $L_0(G^*)$ by letting the $x_e$'s assume all possible
values in $\fq$. Similarly we may replace $y$ by $-y$. From a
particular rank $r$ evaluation of $L_0(G)$ over $\fq$ we can apply the
reasoning in the first proof of Theorem~\ref{thm:kn} to get the
recurrence (\ref{eq:gast}). $\ \Box$

Theorem~\ref{thm:genmcw} provides a simple recursive procedure for
computing $g_G(q)$ for a graph with ``few'' missing edges (and hence
many apexes). For instance, for $n>k$ let $K_n-K_k$ denote the
complete graph $K_n$ on the vertex set $[n]=\{1,2\dots,n\}$ with all
edges $ij$ removed where $i,j\in\{1,2,\dots,k\}$, $i\neq j$. When
$k=n-1$, $L_0(K_n-K_k)$ is just a generic diagonal matrix, so we get
  $$ h(K_n-K_{n-1},r)={n-1\choose r}(q-1)^r. $$
Hence in principle Theorem~\ref{thm:genmcw} can be used iteratively to
compute $h(K_n~-~K_k,r)$ for any $n,k,r$. In particular, it follows that
$g_{K_n-K_k}(q)\in\mathbb{Z}[q]$, verifying Kontsevich's conjecture in
this case. When $k=1$ we have $K_n-K_1=K_n$, which was dealt with in
Theorem~\ref{thm:kn}. When $k=2$ we have $K_n-K_2=K_n-e$, the complete
graph $K_n$ with one edge removed. We compute $g_{K_n-e}(q)$ by
another method in Theorem~\ref{thm:mstar} (the case $s=1$). For $3\leq
k\leq 5$ we have used Theorem~\ref{thm:genmcw} to produce the
following conjecture. One could easily extend this conjecture to
other small values of $k$, but what would be more interesting is a  
conjecture for general $n$ and $k$. 

\begin{conjecture}
\emph{We have}
  \beas 
g_{K_{2m}-K_3}(q) & = & q^{m(m-1)}(q-1)(q^3-1)\cdots
  (q^{2m-5}-1)\\[.05in] & & 
  \quad \cdot (q^{4m-3}-4q^{2m-4}+3q^{2m-5}-q^{2m-6}+1),\ m\geq 2  
  \\[.2in] 
g_{K_{2m+1}-K_3}(q) & = & q^{m^2+m-3}(q-1)(q^3-1)\cdots
  (q^{2m-3}-1)\\[.05in] & & 
  \quad \cdot (q^{2m-1}-3q+2),\ m\geq 2\\[.2in] 
g_{K_{2m}-K_4}(q) & = & q^{m(m-1)}(q-1)(q^3-1)\cdots
  (q^{2m-5}-1)\\[.05in] & & \ 
   \cdot (q^{4m-10}-7q^{2m-6}+8q^{2m-7}-3q^{2m-8}+1),\ m\geq 3\\[.2in]
g_{K_{2m+1}-K_4}(q) & = & q^{m^2+m-4}(q-1)(q^3-1)\cdots
  (q^{2m-5}-1)\\[.05in] & & 
  \quad \cdot
  (q^{4m-6}-8q^{2m-3}+9q^{2m-4}-4q^{2m-5}+q^{2m-6}+4q-3),\\
  & & \qquad\qquad m\geq 2\\[.2in]
g_{K_{2m}-K_5}(q) & = & q^{m(m-1)}(q-1)(q^3-1)\cdots
  (q^{2m-7}-1)\\[.05in] & & \ \
   \cdot (q^{6m-19}-16q^{4m-14}+25q^{4m-15}
  -16q^{4m-16}+5q^{4m-17}-q^{4m-18}\\[.05in] & & \qquad\qquad
    +q^{2m-6}+11q^{2m-8}-15q^{2m-9}+6q^{2m-10}-1),\ m\geq 3\\[.2in]
g_{K_{2m+1}-K_5}(q) & = & q^{m^2+m-5}(q-1)(q^3-1)\cdots
  (q^{2m-5}-1)\\[.05in] & & 
    \quad \cdot
  (q^{4m-9}-15q^{2m-5}+24q^{2m-6}-15q^{2m-7}+4q^{2m-8}\\
  & & \qquad\qquad +5q-4),\ m\geq 3. 
\eeas 
\end{conjecture}

To prove this conjecture for a particular value of $k$, one could try
to guess a formula for $h(K_n-K_k,r)$ and then verify that it
satisfies the recurrence (\ref{eq:gast}) (with appropriate initial
conditions). 

For our second generalization of Theorem~\ref{thm:kn}, we need to
consider the inequivalent nondegenerate symmetric scalar products on
the space $\fq^n$. Standard results in orthogonal geometry over $\fq$
(essentially equivalent to the results used in the second proof of
Theorem~\ref{thm:kn}) show that there are two such scalar products
when $q$ is odd. They are defined as follows, where
$a=(a_1,\dots,a_n)$, $b=(b_1,\dots,b_n)\in\fq^n$.
  \beq \langle a,b\rangle_+  =  ab^t  =  \sum
    a_ib_i \label{eq:adotb} \eeq 
  $$ \langle a,b\rangle_-  =  aDb^t 
    =  \alpha a_1b_1 + \sum_{i=2}^n a_ib_i. $$
Here $D$ and $\alpha$ have the same meaning as in the second proof of
Theorem~\ref{thm:kn}, so in particular $\alpha$ is a nonsquare in
$\fq$.

Similarly when $q$ is even, if $n$ is odd then all nondegenerate
symmetric scalar products are equivalent to (\ref{eq:adotb}). When
$n=2m$ we have (\ref{eq:adotb}) together with
   $$ \langle a,b\rangle_- = aEb^t = \sum_{i=1}^m (a_{2i-1}b_{2i}
     +a_{2i}b_{2i-1}), $$
where $E$ is defined in the second proof of Theorem~\ref{thm:kn}. 

Now suppose that $G$ is a graph with vertex set $\{v_1,\dots,v_{n+1}\}$
such that $v_{n+1}$ is an apex. Let $b_G^+(q)$ (respectively,
$b_G^-(q)$) denote the number of ordered bases $(u_1,\dots,u_n)$ of
$\fq^n$ such that
  $$ \langle u_i,u_j\rangle_+  = 0\ \mbox{(respectively,
      $\langle u_i,u_j\rangle_-=0$)}, $$
whenever $i\neq j$ and $ij\not\in E(G)$. Such an ordered basis forms
the rows of a matrix $A\in\gln$ such that $(AHA^t)_{ij}=0$ whenever
$i\neq j$ and $ij\not\in E(G)$, where $H=I,D$, or $E$ depending on
which of the scalar products we are considering. It follows that the
second proof of Theorem~\ref{thm:kn} extends \emph{mutatis mutandis}
to give the following result.

\begin{theorem} \label{thm:orbas}
\emph{Let $G$ be as above. If $q$ is odd or if $q$ is even and $n$ is
even, then}
  $$ g_G(q) = \frac{b_G^+(q)}{\#\Omega^+(n,q)} +
  \frac{b_G^-(q)}{\#\Omega^-(n,q)}. $$
\emph{If $q$ is even and $n$ is odd, then}
  $$ g_G(q) = \frac{b_G^+(q)}{\#\Omega(n,q)}. $$
\end{theorem}

As an example of the use of Theorem~\ref{thm:orbas},
let $K_{1,s}$ denote the \emph{star} consisting of one vertex
connected to $s$ other vertices, and let $G=K_{n+1}-K_{1,s}$ for
$n>s+1$. In other words, $G$ consists of $K_{n+1}$ with $s$ edges
removed which are incident to a common vertex. In particular,
$K_{n+1}-K_{1,1} = K_{n+1}-K_2$, the special case $k=2$ of $K_n-K_k$
considered above (with $n$ replaced by $n+1$).

\begin{theorem} \label{thm:mstar}
\emph{We have}
  \bea g_{K_{2m-1}-K_{1,s}}(q) & = & q^{m^2+m-s-1}(q-1)(q^3-1)\cdots
       (q^{2m-3}-1)\nonumber\\ & & \quad\cdot(q^{2m}-q^s-q+1),\quad
  s\leq  2m-3\label{eq:kis1}\\[.1in] 
  g_{K_{2m}-K_{1,s}}(q) & = & q^{m(m-1)}(q-1)(q^3-1)\cdots
   (q^{2m-3}-1)\nonumber \\ & & \quad \cdot (q^{2m-1-s}-1),\quad s\leq
  2m-2. \label{eq:kis2} \eea 
\end{theorem}

\textbf{Proof.} According to Theorem~\ref{thm:orbas} we need to count
the number of ordered bases $(u_1,\dots,u_n)$ of $\fq^n$ satisfying
  $$ \langle u_1,u_2\rangle_+ =\cdots=\langle u_1,u_{s+1}\rangle_+=0,
  $$ 
as well as the number of ordered bases $(u_1,\dots,u_n)$ of $\fq^n$
satisfying 
  $$ \langle u_1,u_2\rangle_- =\cdots=\langle u_1,u_{s+1}\rangle_-=0
  $$ 
(except that when $q$ is even and $n$ is odd we only have one type of
scalar product). Let $u_1^\perp$ denote the set of all vectors
orthogonal to $u_1$, with respect to whatever scalar product is under
consideration. We always have $\dim u_1^\perp=n-1$. Once we have
chosen $u_1$, if $u_1\not\in
u_1^\perp$ then there are $(q^{n-1}~-~1)(q^{n-1}~-~q)\cdots
(q^{n-1}~-~q^{s-1})$ 
choices for $u_2,\dots,u_{s+1}$, and then
$(q^n~-~q^{s+1})\cdots(q^n~-~q^{n-1})$ choices for
$u_{s+2},\dots,u_n$. On the other hand, if 
$u_1\in u_1^\perp$ then there are $(q^{n-1}-q)(q^{n-1}-q^2)\cdots
(q^{n-1} -q^s)$ choices for $u_2,\dots,u_{s+1}$, and then
$(q^n-q^{s+1})\cdots (q^n-q^{n-1})$ choices for $u_{s+2},\dots, u_n$
as before. Hence to complete the computation we need to know the
number $N_{\pm}(n)$ of $u_1$ for which $u_1\in u_1^\perp$, i.e,
$\langle u_1,u_1 \rangle_\pm= 0$. When $q$ is even it is easy to
compute $N_\pm(n)$, while when $n$ is odd this number appears e.g.\ in
\cite[Thms.\ 65 and 66]{dickson}\cite[Thms.\
1.26 and 1.37]{wan}. The values are
  $$ \begin{array}{rcll} N_+(n) = N(n) & = & q^{n-1}, & q\
  \mathrm{odd\ and}\ n\ \mathrm{odd}\\[.05in] 
  N_+(n) & = & q^{n-1}+q^{\frac n2}-q^{\frac n2 -1}, & \mbox{either
  $q$ is odd and $n \equiv 0$ (mod 4),}\\ & & &\ \ \mbox{or $q\equiv 1$
  (mod 4) and 
  $n\equiv 2$ (mod 4)}\\[.05in]
  N_+(n) & = & q^{n-1}-q^{\frac n2}+q^{\frac n2 -1}, & \mbox{either
  $q$ is odd and $n \equiv 2$ (mod 4),}\\ & & &\ \ \mbox{or $q\equiv 1$
  (mod 4) and 
  $n\equiv 0$ (mod 4)}\\[.05in]
  N_-(n) & = & q^{n-1}-q^{\frac n2}+q^{\frac n2 -1}, & \mbox{either
  $q$ is odd and $n \equiv 0$ (mod 4),}\\ & & &\ \ \mbox{or $q\equiv 1$
  (mod 4) and 
  $n\equiv 2$ (mod 4)}\\[.05in]
  N_-(n) & = & q^{n-1}+q^{\frac n2}-q^{\frac n2 -1}, & \mbox{either
  $q$ is odd and $n \equiv 2$ (mod 4),}\\ & & &\ \ \mbox{or $q\equiv 1$
  (mod 4) and 
  $n\equiv 0$ (mod 4)}\\[.05in]
  N(n) & = & q^{n-1}, & \mbox{$q$ even and $n$ odd}\\[.05in]
  N_+(n) & = & q^{n-1}, & \mbox{$q$ even and $n$ even}\\[.05in]
  N_-(n) & = & q^n, & \mbox{$q$ even and $n$ even}.
  \end{array} $$
It is now a routine computation (which we omit) to obtain the stated
formulas (\ref{eq:kis1}) and (\ref{eq:kis2}). $\ \Box$

\section{Some related negative results.} \label{sec:neg}
In Theorem~\ref{thm:apex} we showed that the Kontsevich conjecture for
graphs with apexes is equivalent to counting nonsingular symmetric
matrices over $\fq$ with specified ``holes'' (entries equal to 0),
with no holes on the main diagonal. A related problem that comes to
mind is the case of arbitrary matrices, rather than symmetric
matrices. Thus let $S$ be any subset of $[n] \times [n]$, and let
$h_S(q)$ denote the number of matrices $A\in\gln$ whose support is
contained in $S$, i.e., $A_{ij}=0$ whenever $(i,j)\not\in S$. For
instance, 
  $$ h_{[n]\times [n]}(q) = \#\gln =
  (q^n-1)(q^n-q)\cdots(q^n-q^{n-1}). $$ 

\textbf{Question.} Is the function $h_S(q)$ always a polynomial in
$q$? 

According to Kontsevich (private communication), a negative answer
follows from the existence of the Fano plane $F$ (the projective plane
of order two, with three points on a line and seven points in all). We
have not been able to understand this remark of Kontsevich. However,
if we take $n=7$ and let $S\subset [7]\times[7]$ be the support of
the incidence matrix of $F$ (so $\#S=21$), then Stembridge has
shown that $h_S(q)$ is \emph{not} a polynomial. More precisely,
  $$ h_S(q) = \left\{ \begin{array}{ll}
      q^{21}-q^{20}-q^{19}-14q^{18}-7q^{17}+176q^{16}
      +8q^{15}-1860q^{14}\\ \ \ +5603q^{13}-8880q^{12}+9010q^{11}
      -6110q^{10}+2603q^9\\ \ \ -428q^8-248q^7+208q^6-72q^5+13q^4-q^3, &
     q\ \mathrm{odd}\\[.1in] 
      q^{21}-q^{20}-q^{19}-14q^{18}-7q^{17}+175q^{16}+21q^{15}
    -1938q^{14}\\ \ \ +5889q^{13}-9595q^{12}+10297q^{11}-7826q^{10}
      +4319q^9\\ \ \ -1715q^8+467q^7-78q^6+6q^5, & q\ \mathrm{even}.
     \end{array} \right. $$ 
Moreover, Stembridge has also verified that $S$ is the \emph{smallest}
counterexample to the polynomiality of $h_T(q)$, in the sense that
$h_T(q)$ is a polynomial whenever $\#T \leq 21$ or whenever $n\leq 7$,
except when $T$ can be transformed to $S$ by row and column
permutations. 


Now let $S$ be a \emph{symmetric} subset of $[n]\times [n]$, i.e.,
$(i,j)\in S\Leftrightarrow (j,i)\in S$. Define $k_S(q)$ to be the
number of invertible $n\times n$ \emph{symmetric} matrices over
$\mathbb{F}_q$ whose support is contained in $S$. Suppose that $T$ is
a subset of $[n]\times [n]$ for which $h_T(q)$ is not a 
polynomial. (The discussion above shows that we can take $n=7$.) Let 
$A$ be an $n\times n$ matrix over $\fq$.  Then $A$ is counted by
$h_T(q)$ if and only if the $2n\times 2n$ matrix
  $$ B = \left[ \begin{array}{cc} 0 & A\\ A^t & 0 \end{array}
   \right] $$
is a nonsingular \emph{symmetric} matrix with support contained in
  $$ S = \{ (i,j+n), (j+n,i)\,:\,(i,j)\in T\}\subset [2n]\times
  [2n]. $$ 
Hence $k_S(q)$ is not a polynomial.
Unfortunately all the main diagonal elements
of examples of the form $\left[ \begin{array}{cc} 0 & A\\ A^t & 0
\end{array} \right]$ are holes, so we cannot use
Theorem~\ref{thm:apex} to deduce that we have a counterexample to 
Kontsevich's conjecture. 

As pointed out by Stembridge, there are even simpler examples of
symmetric sets $S\subset [n]\times[n]$ for which $k_S(q)$ is not a
polynomial. Any symmetric matrix of odd order $n$ with 0's on the main
diagonal over a field of characteristic 2 is singular. Hence for $n$
odd we can choose $S$ to be any subset of $[n]\times [n]$ that
includes no element of the form $(i,i)$ and that contains at least one
transversal (i.e., a subset $(i,w(i))$ where $w$ is a permutation of
$[n]$). (This last condition is equivalent to $k_S(q)\neq 0$ for some
$q$.) Then $k_S(q)=0$ for $q=2^m$, but $k_S(q)\neq 0$ for some $q$, so
$k_S(q)$ is not a polynomial. Since the prime 2 plays such a special
role in this example, perhaps the function $h_S(q)$ or $k_S(q)$ is a
polynomial in $q$ for \emph{odd} $q$. A good place to look for a
counterexample to this suggestion would be when $S$ is the support of
a projective plane of odd order, but even for the plane of order 3 we
are unable to compute $h_S(q)$.

There are various natural generalizations of Kontsevich's
conjecture. For instance, the spanning trees of a connected graph $G$
form the bases of the graphic matroid associated with $G$ (see e.g.\
\cite[1.3.B]{crapo}). Thus if $M$ is any matroid on the set
$\{e_1,\dots,e_s\}$, 
then define in complete analogy to (\ref{eq:qgdef})
  $$ Q_M(x) = \sum_B x^B, $$
where $B$ ranges over all bases of $M$. Let $g_M(q)$ denote the number
of ways of substituting elements of $\fq$ for the variables of $Q_M(x)$
such that $Q_M(x)\neq 0$. We can generalize Kontsevich's conjecture by
asking whether $g_M(q)$ is always a polynomial function of
$q$. There are, however, very simple counterexamples. For instance, if
$M$ is the four-point line so that
  $$ Q_M(x) = x_1x_2+x_1x_3+x_2x_3+x_1x_4+x_2x_4+x_3x_4, $$
then it can be shown that 
  $$ g_M(q) = \left\{ \begin{array}{rl} 
     q(q-1)(q^2-1), & q\equiv 1\ (\mathrm{mod}\ 3)\\[.1in]  
     q(q-1)(q^2+1), & q\equiv 2\ (\mathrm{mod}\ 3)\\[.1in]  
     q^3(q-1), & q\equiv 0\ (\mathrm{mod}\ 3). \end{array} \right. $$
Matroid theorists will notice that the four-point line is not a
regular (or unimodular) matroid, but every graphic matroid is
regular. Hence it is natural to ask whether $g_M(q)$ might be a
polynomial for regular matroids $M$. However, Stembridge has shown
that for the regular matroid $M$ called R10 in Oxley's book
\cite{oxley}, $g_M(q)$ is not a polynomial.

Although for the four-point line $g_M(q)$ is not a polynomial, note
that it is a \emph{quasipolynomial}, i.e., for some $N>0$ (here $N=3$)
it is a polynomial on the different residue classes modulo $N$. Thus
it might be interesting to consider for which matroids $M$ (or for
even more general varieties than the zeros of $Q_M(x)$) is $g_M(q)$ a
quasipolynomial. In particular, if Kontsevich's conjecture is false,
is it at least true that $g_G(q)$ is a quasipolynomial?

\textsc{Acknowledgment.} I am grateful to Jeff Lagarias for calling my
attention to the conjecture of Kontsevich and for providing me with a
copy of his notes of Kontsevich's lecture and a related lecture of I.
M. Gelfand. I also wish to thank Tim Chow and John Stembridge for some
useful discussions.

\end{document}